\documentclass[12pt]{article}
\usepackage{amsmath}
\usepackage{amssymb}
\usepackage{amsfonts}

\newtheorem{thm}{Theorem}[section]
\newtheorem{pro}[thm]{Proposition}
\newtheorem{cor}[thm]{Corollary}
\newtheorem{lem}[thm]{Lemma}
\newtheorem{dfn}[thm]{Definition}
\newtheorem{fct}[thm]{Fact}

\newtheorem{rmk}[thm]{Remark}
\newtheorem{ass}[thm]{Assumption}

=22cm\headheight=0cm\topskip=0cm\topmargin=0cm

\def\cal{\mathcal {}}
\def\qed{\hfill$\Box$}

\newcommand{\M}{\sf M}

\begin{document}
\def\dis{\displaystyle}

\begin{center}
{{\sc 
correspondences between model theory and banach space theory  }}
\vspace{10mm}

{
{\bf  Karim Khanaki}} \vspace{3mm}

{\footnotesize
  Faculty of Fundamental Sciences, Arak University of Technology,
 \\
P.O. Box 38135-1177, Arak, Iran; and
\\ School of Mathematics,
Institute for Research in Fundamental Sciences (IPM), \\ P.O. Box
19395-5746, Tehran, Iran; \\
 e-mail:
khanaki@arakut.ac.ir} \vspace{5mm}
\end{center}

{\sc Abstract.}
{\small In \cite{K3} we pointed out the correspondence between a
result of Shelah in model theory, i.e. a theory is unstable if and
only if it has IP or SOP, and the well known compactness theorem
of Eberlein and \v{S}mulian in functional analysis. In this
paper, we relate a {\em natural} Banach space $V$ to a formula
$\phi(x,y)$, and show that $\phi$ is stable (resp NIP, NSOP) if
and only if $V$ is reflexive (resp Rosenthal, weakly sequentially
complete) Banach space. Also,
 we present a proof of the Eberlein-\v{S}mulian theorem by a
model theoretic approach using Ramsey theorems which is
illustrative to show some correspondences  between model theory
and Banach space theory. }

\medskip



\tableofcontents

\section{Introduction}
 In \cite{Sh} Shelah introduced the independence
property (IP) for 0-1 valued formulas and defined the strict
order property as complementary to the independence property: a
theory is unstable iff it has IP or SOP. On the other hand, a
well known fact in functional analysis, the Eberlein-\v{S}mulian
theorem, states that a subset of a Banach space is not  weakly
precompact iff it has a sequence without any weak Cauchy
subsequence or it has a weak Cauchy sequence with no weak limit.
In fact, Shelah's result corresponds to
 the Eberlein-\v{S}mulian theorem. This was noticed in \cite{K3} and some various forms of
definability of types for NIP models were proved.

The relation between model theory and Banach space theory is
rather deep (e.g. see \cite{BBHU} and its references). So it
would be desirable to have clearly understood channels of
communication between technical complexity that both fields have
attained in the last thirty years, so that techniques from one
field might become useful in the other.

In this paper we purpose to give a proof of the
Eberlein-\v{S}mulian theorem with a model theoretic point of view
and investigate correspondence between model theory and Banach
space theory. We relate a {\em natural} Banach space $V$ to a
formula $\phi(x,y)$, and show that $\phi$ is stable (resp NIP,
NSOP) if and only if $V$ is reflexive (resp Rosenthal, weakly
sequentially complete) Banach space.

However, we believe that the main goal of this paper, if it is
achieved, is to show that the correlation  between model theory
and Banach space theory is much more than what is known so far,
and there are many connections between model theoretic
classification and Banach space classification which can be
studied in future works.

It is worth recalling another lines of research. In \cite{Iba} and
\cite{S} the relationship between NIP and Rosenthal's dichotomy
were noticed in the contexts of $\aleph_0$-categorical structures
in continuous logic and classical first order setting,
respectively. The relationship between NIP in integral logic and
Talagrand's stability was studied in  \cite{K}. The above
correspondence was noticed in \cite{K3} and some various forms of
definability of types for NIP models were proved.

This paper is organized as follows: In the next section we
present a proof of the Eberlein-\v{S}mulian theorem and study its
connection to the theorem of definability of types.
In the third section, we
investigate some correspondences between model theory and Banach
space theory.

\bigskip\noindent
{\bf Acknowledgements.} I would like to thank the Institute for
Basic Sciences (IPM), Tehran, Iran. Research partially supported
by IPM grants nos 92030032 and 93030032.

\noindent\hrulefill

\section{Eberlain-\v{S}mulian and Ramsey theorems} \label{Eberlein}
The Eberlein-\v{S}mulian theorem states that the following three
versions of precompactness are equivalent for the weak topology of
a Banach space. Let $A$ be a subset of a topological space $X$,
then
\begin{itemize}
             \item [(i)] The set $A$ is {\em precompact} if its closure is compact.
             \item [(ii)] The set $A$ is {\em sequentially precompact} if each sequence of elements of $A$ has a subsequence converging to an element of $X$.
              \item [(iii)] The set $A$ is {\em countably precompact} if each sequence of elements of $A$ has a cluster point in $X$.
\end{itemize}

To prove this theorem, first we give a proof for the Banach space
$C(X)$ which $X$ is a compact topological space. Then the general
case follows easily from this case (see Remark~\ref{rem} below).
It is clear that (ii) implies (iii). We prove that (i) implies
(ii) using a stronger version of Ramsey's theorem (see
Theorem~\ref{Ramsey} below). The proof of (iii)~$\Rightarrow$~(i)
uses the theorem of definability of types and a combinatorial
lemma.

\subsection{Pointwise convergence}
First we review some notions and results for the topology of
pointwise convergence. If $X$ is any set and $A$ a subset of
$\mathbb{R}^X$, then the topology of {\em pointwise convergence}
on $A$ is that inherited from the usual product topology of
$\mathbb{R}^X$; that is, the coarsest topology on $A$ for which
the map $f\mapsto f(x):A\to \mathbb R$ is continuous for every
$x\in X$. A typical neighborhood of a function $f$ is determined
by a finite subset $\{x_1,\ldots, x_n\}$ of $X$ and $\epsilon>0$
as follows:  $$U_f(x_1,\ldots,x_n;\epsilon)=\{g\in\Bbb R^X:
|f(x_i)-g(x_i)|<\epsilon\mbox{ for }i\leqslant n\}.$$

\begin{ass} In this paper (countable, sequential) precompactness
means (countable, sequential) precompactness with respect to the
topology of pointwise convergence.  Also, we will say that a
sequence or net is  convergent  if it is convergent for the
topology of pointwise convergence. Otherwise, we explicitly state
that what is our desired topology; e.g. weak, weak* or norm
topology.
\end{ass}

Now, recall the following standard fact from functional analysis.
In fact, it is a topological presentation of stability in model
theory (see \cite{K3}). (See the appendix for its proof.)

\begin{fct}[Grothendieck's criterion] \label{Criterion} 
 Let $X$ be a compact topological space. Then the following are equivalent for a norm-bounded subset
$A\subseteq C(X)$:
\begin{itemize}
             \item [(i)]  The set $A$ is precompact in $C(X)$.
             \item [(ii)] For every sequences $\{f_n\}\subseteq A$ and
$\{x_n\}\subseteq X$, we have $$\lim_n \lim_m f_n(x_m) =\lim_m
\lim_n f_n(x_m),$$  whenever both limits exist.
\end{itemize}
\end{fct}

Here, we will present a form of Ramsey's theorem which we use in
this paper. For this, we give some notations. Let $[\mathbb{N}]$
denote all infinite subsequences of $\mathbb N$. For $n\in\mathbb
N$, $[\mathbb{N}]^n$ denotes all finite subsequences of $\mathbb
N$ of length $n$. If $M\in[\mathbb N]$ we use similar notation,
$[M]$ and $[M]^n$ to denote all subsequences (or all length $n$
subsequences) of $M$. Suppose we can color, using colors $R$ and
$B$, all infinite subsequences of $\mathbb N$. Thus $M\in
[\mathbb N]$ implies $M\in R$ or $M\in B$. The following Ramsey's
theorem shows that if the set $R$ has a ``good" property then
there exists $M\in [\mathbb N]$ so that either for all $N\in [M]$,
$N\in R$ or for all $N\in [M]$, $N\in B$. To be more precise we
topologize $[\mathbb N]$ by the product topology; thus a basic
open set in $[\mathbb N]$ is of the form
$$O(n_1,\ldots,n_k)=\{M=(m_i)\in [\mathbb N] : m_i = n_i \mbox{
for } i\leqslant k\}$$
 where $n_1<\cdots <n_k$ is arbitrary.

\begin{thm}[Ramsey theorem]  \label{Ramsey} Assume that $\mathcal A\subseteq [\mathbb N]$ is Borel. Then there
is an $N\in [\mathbb N]$ so that either $[N]\subseteq \mathcal A$
or $[N]\cap \mathcal A=\emptyset$.
\end{thm}

The above theorem is due to Galvin and Prikry in \cite{GP}.

\begin{dfn} \label{indep}  A sequence $\{f_n\}$ of real valued functions on a set $X$ is said to be independent if there
exist real numbers $s<r$ such that
$$\bigcap_{n\in P} f_n^{-1}(-\infty,s)\cap\bigcap_{n\in M} f_n^{-1}(r,\infty)\neq\emptyset$$
for all finite disjoint subsets $P,M$ of $\mathbb N$.
\end{dfn}

\begin{lem}[Rosenthal] \label{lemma-1} Let $X$ be a compact space and $F\subseteq C(X)$ a bounded subset. Then the following conditions are equivalent:
\begin{itemize}
             \item [(i)]  $F$ does not contain an independent subsequence.
             \item [(ii)] Each sequence in $F$ has a  convergent subsequence
in $\mathbb R^X$.
\end{itemize}
\end{lem}
{\bf Proof.} (i)~$\Rightarrow$~(ii): Suppose that
$\{f_n\}\subseteq F$ has no convergent subsequence. Therefore,
there are rational numbers $r>s$ such that for all infinite
subset $M\subseteq \mathbb N$ there exists $x\in X$ so that $x$
belongs to infinitely many $A_n=f_n^{-1}(-\infty,s)$'s, $n\in M$
and to infinitely many $B_n=f_n^{-1}(r,\infty)$'s, $n\in M$.
(Indeed, if this were not true, let $(r_i,s_i)_{i=1}^\infty$
(with $r_i>s_i$) be dense in $\mathbb R^2$ and inductively choose
$M_{i+1}\in[M_i]$ so that the above conditions fail. If $M$ is a
diagonal of $M_i$'s then $\{f_n\}_{n\in M}$ is convergent.) Now,
let
$${\mathcal A}=\big\{L=(x_i)_{i=1}^\infty\in X^{\mathbb N}:
\bigcap_{i=1}^k A_{l_{2i-1}}\cap\bigcap_{i=1}^k
B_{l_{2i}}\neq\emptyset\mbox{ for all }k\in\mathbb{N}\big\}.$$
Then $\mathcal A$ is closed, because $\mathcal
A=\bigcap_{k=1}^\infty\mathcal A_k$ where $\mathcal
A_k=\{L=(x_i)_{i=1}^\infty\in X^{\mathbb N}: \bigcap_{i=1}^k
A_{l_{2i-1}}\cap\bigcap_{i=1}^k B_{l_{2i}}\neq\emptyset\}$. (We
note that $\mathcal A_k$'s are closed.) By Ramsey
theorem~\ref{Ramsey} and the above observation, there is some
infinite subset $L\subseteq\mathbb N$ such that $[L]\subseteq
\mathcal A$. Clearly $\{f_{l_{2i}}\}_{i<\omega}$ is independent.

(ii)~$\Rightarrow$~(i) is an easy exercise. (Indeed, if $\{f_n\}$
was independent and $r>s$ were witness, let $\{f_{n_k}\}$ be a
convergent subsequence of it, and let $M\subseteq\omega$ be
infinite and coinfinite in $\{n_k:k<\omega\}$. Let $M_n$ be the
initial segment of $M$ where $|M_n|=n$, and let $x_M$ be a cluster
point of $\{x_{M_n}:n<\omega\}$ where for all $n$,
$f_{n_k}(x_{M_n})>r$ if $k\in M_n$ and $f_{n_k}(x_{M_n})<s$ if
$k\in\{1,\ldots,\max(M_n)\}\setminus M_n$. Then $f_{n_k}(x_M)$
converges to two different numbers, a contradiction.) \qed

\medskip

Haskell P. Rosenthal \cite{Ros} used the above lemma for proving
his famous $\ell_1$ theorem: a sequence in a Banach space is
either `good' (it has a subsequence which is weakly Cauchy) or
`bad' (it contains an isomorphic copy of $\ell_1$). we will
shortly  discuss this topic.

\medskip
we continue the discussion of some the topics The
following lemma is a generalization of a model theoretic fact,
i.e. IP implies OP.

\begin{lem}[IP~$\Rightarrow$~OP] \label{IP->OP}
Let $X$ be compact and $A\subseteq C(X)$ be  precompact (in
$C(X)$). Then every sequence in $A$ has a convergent subsequence
in $\mathbb{R}^X$, and so in $C(X)$.
\end{lem}
{\bf Proof.} Suppose, if possible, that the sequence
$\{f_n\}\subseteq A$ has no convergent subsequence. Thus, by
Lemma~\ref{lemma-1},  $\{f_n\}$ is independent, i.e. there are
$r>s$ such that for all finite disjoint $P,M\subseteq\mathbb{N}$,
we have
 $$\Big\{x\in X:\big( \bigwedge_{i\in P}f_i(x)\leqslant s\big)\wedge\big(\bigwedge_{i\in M}f_i(x)\geqslant
r\big)\Big\}\neq\emptyset.$$ Since $X$ is compact, in the
definition of independent sequence, one can assume that
$P\subseteq \mathbb N$ is infinite and $M=\mathbb N\setminus P$.
 Now, a straightforward adaptation of a classical result in
 model theory, i.e. IP implies OP,  shows that there are subsequences of $\{f_n\}$ (still
 denoted by $\{f_n\}$) and $\{x_m\}$ in $X$ such that
 $\lim_m\lim_nf_n(x_m)\geqslant r>s\geqslant\lim_n\lim_mf_n(x_m)$.
(Indeed, suppose that $(x_P)_{P\in 2^\omega}$ witness IP. Given $i
<\omega$, let $P_i:\omega\to 2$ such that $P_i(j) = 0$ if and only
if $i\leqslant j$. Then we have $f_i(x_{P_j})\leqslant s$ iff
$P_j(i)=1$ iff $i<j$, and $f_i(x_{P_j})\geqslant r$  iff
$P_j(i)=0$  iff $i\geqslant j$. Take $x_m=x_{p_m}$.)
 This is a contradiction by Fact~\ref{Criterion}. Moreover, since $A$ is
 precompact, the limit of every convergent sequence is
 continuous.
\qed

\medskip

Now, we are ready to give a proof of the Eberlain-\v{S}mulian for
the topology of pointwise convergence on $C(X)$.

\begin{thm}[Eberlain-\v{S}mulian for $C_p(X)$] \label{ES-pointwise} Suppose that $C(X)$ is the Banach
space of all continuous real-valued functions of on a compact
space $X$ and $A\subseteq C(X)$ is norm-bounded. Then for the
topology of pointwise convergence the following are equivalent.
\begin{itemize}
             \item [(i)]  The set $A$ is precompact.
             \item [(ii)] The set $A$ is sequentially precompact.
             \item [(iii)] The set $A$ is countably precompact.
\end{itemize}
\end{thm}
{\bf Proof.} (i)~$\Rightarrow$~(ii) is just Lemma~\ref{IP->OP}.

(iii)~$\Rightarrow$~(i):  Suppose that  $A\subseteq C(X)$ is
countably precompact and norm-bounded.
 Suppose that $f_n\in A$ and $x_n\in X$ form two sequences and
the limits $\lim_n \lim_m f_n(x_m)$ and $\lim_m \lim_n f_n(x_m)$
exist. Let $f$ in $C(X)$ and $x$ in $X$ be cluster points of
$\{f_n\}$ and $\{x_m\}$. Thus,
\begin{align*} \lim_n\lim_m f_n(x_m)  =
\lim_nf_{n}(x) &= f(x)=\lim_mf(x_m)=\lim_m\lim_nf_n(x_m).
\end{align*}
By Fact~\ref{Criterion}, $A$ is precompact. \qed

\medskip

This is the end of the story if one replaces the `topology of
pointwise convergence'  with the `weak topology'. But, it needs
some works:

\subsection{Weak topology} Recall that for a normed space $U$, the
topology generated by $U^*$ is known as the weak topology on $U$.
The next remark states that every normed space with its weak
topology lives inside a space of the form $C(X)$, with the
topology of pointwise convergence, where $X$ is a compact space.

\begin{rmk} \label{rem} {\em For an arbitrary normed space $U$,
write $X$ for the unit ball of the dual space $U^*$, with its
weak* topology. Then $X$ is compact by Alaoglu's theorem and the
natural map $u\mapsto \hat{u}:U\to\mathbb{R}^X$, defined by
setting $\hat{u}(x)=x(u)$ for $x\in X$ and $u\in U$, is a
homeomorphism between $U$, with its weak topology, and its image
$\widehat{U}$ in $C(X)$, with the topology of pointwise
convergence. }
\end{rmk}

If we show that the direction (i)~$\Rightarrow$~(ii) in
Theorem~\ref{ES-pointwise} holds for any subspace $Y$ of $C(X)$,
i.e. every precompact subset of subspace $Y$ is sequentially
precompact in $Y$, then by Remark~\ref{rem}, for every Banach
space (or even normed space), the weak precompactness implies
weak sequential precompactness. Indeed, if $A\subseteq Y$ be
 precompact in $Y$, then $A$ is countably
 precompact in $Y$, and so is countably
 precompact in $C(X)$. Therefore, if $\{f_n\}\subseteq
A$, by the direction (i)~$\Rightarrow$~(ii) for $C(X)$, there
exists $\{g_n\}\subseteq\{f_n\}$ and $f\in C(X)$ such that $g_n\to
f$. Since $A$ is countably   precompact in $Y$, the cluster point
of $\{g_n\}$, i.e. $f$, is in $Y$. Thus the direction
(i)~$\Rightarrow$~(ii) holds for any subspace $Y\subseteq C(X)$.

In the next section we present a proof of the direction
(iii)~$\Rightarrow$~(i) using the theorem of definability of types
and a well known combinatorial result.


\section{Definability of types}
 For  convex subsets, one can show that the direction (iii)~$\Rightarrow$~(i)
is a consequence of  a well known fact in model theory, that the
theorem of definability of types. This theorem says that for a
formula $\phi(x,y)$ stable in a model $M$, and every type $p\in
S_\phi(M)$ there is a sequence $\psi_n(y)$ of the convex
combinations of $\phi(a,y)$'s, $a\in M$, such that $\psi_n(y)$ is
uniformly convergent to a (continuous) function $\psi(y)$ where
$\phi(x,b)^p=\psi(b)$ for all $b\in M$. The following proof is a
straightforward translation of the proof of definability of types
for continuous logic, as can be found in \cite[Appendix~B]{BU}:

\begin{fct}[Definability of types] \label{Definability}
Let $X$ be a compact space and $A\subseteq C(X)$ norm-bounded and
countably  precompact. Then every point of the closure of $A$ is
a uniform limit of a sequence in the convex hull of $A$, denoted
by $\textrm{conv}(A)$.
\end{fct}
{\bf Proof.} Since $A$ is countably
 precompact in $C(X)$, it is stable in the sense of model
theory, i.e. the condition (ii) in Fact~\ref{Criterion} holds.
 With out lose of generality we can assume that all functions in $A$ are
$[0,1]$-valued.  Let $f\in\bar{A}$. We claim that for any
$\epsilon>0$, there is a finite sequence $(f_i:i<n_\epsilon)$ in
$A$ such that for all $x,y\in X$, if for all $i<n_\epsilon$,
$|f_i(x)- f_i(y)|\leqslant\epsilon$ then $|f(x)- f(y)|\leqslant
3\epsilon$. If not, by induction on $n$ one can find $f_n\in A$,
$x_n,y_n\in X$ as follows.
  At each step, there are by assumption $x_n,y_n \in X$ such that
  $|f_i(x_n) - f_i(y_n)|\leqslant \epsilon$ for all $i < n$, and yet
  $|f(x_n)- f(y_n)|> 3\epsilon$.
  Once these choices are made, since $f\in\bar{A}$ we may therefore find
  $f_n \in A$ such that $|f_n(x_i)-f_n(y_i)|>3\epsilon$ for all $i \leq n$.
When the construction is complete, we have
$\lim_i\lim_n|f_i(x_n)-f_i(y_n)|\leqslant\epsilon<3\epsilon
\leqslant\lim_n\lim_i|f_i(x_n)-f_i(y_n)|$, a contradiction. Now,
we define increasing function
$g_\epsilon:[0,1]^{n_\epsilon}\to[0,1]$ by
$g(\bar{u})=\sup\{f(x):x\in X \mbox{ and } f_i(x)\leqslant u_i
\mbox{ for all } i<n_\epsilon\}$. Clearly we can find a continuous
increasing function $h_\epsilon$ such that
$g_\epsilon(\bar{u})\leqslant h_\epsilon(\bar{u})\leqslant
g_\epsilon(\bar{u}+\epsilon)$. Thus,
$|f(x)-h_\epsilon(f_i(x):i<n_\epsilon)|\leqslant 3\epsilon$. As
$\epsilon$ is arbitrary, for $\epsilon=\frac{1}{n}$, there is a
continuous function $f_n(x)=h_\epsilon(f_i(x):i<n_\epsilon)$ such
that $|f(x)-f_n(x)|<3\epsilon$ for all $x\in X$, and so $f_n\to
f$ uniformly. \qed

\medskip

\textbf{The direction (iii) $\Rightarrow$ (i) for convex sets.}
Let $Y$ be any subset of $C(X)$. If a {\em convex} set
$A\subseteq Y$ is  countably precompact in $Y$, then $A$ is
  countably precompact in $C(X)$, so $\bar A$, the closure
of $A$ in $C(X)$, is compact  (see Theorem~\ref{ES-pointwise}).
Now if $x\in \bar A$, by the theorem of definability of types
(Fact \ref{Definability}), there is a sequence $(x_n)$ in $A$
converging to $x$; but $(x_n)$  must have a cluster point in $Y$,
and (because the topology is Hausdorff) this cluster point can
only be $x$. Accordingly $\bar{A}\subseteq Y$ and is the closure
of $A$ in $Y$. Thus $A$ is precompact in $Y$. Therefore, by
Remark~\ref{rem} above, a convex subset $A$ (of a normed space
$Y$) is weakly precompact if it is weakly countably precompact.

\subsection{Pt\'{a}k's lemma and stability}
By a combinatorial result due to Pt\'{a}k's, one can show that
the theorem of definability of types implies the direction (iii)
$\Rightarrow$ (i) of the main theorem. First we need some
definitions.

A convex mean on $\Bbb N$ is a function $\mu:\Bbb N\to[0,1]$ such
that (1) $\sum_{i=0}^\infty \mu(i)=1$, (2)
$\textrm{supp}(\mu)=\{i:\mu(i)>0\}$ is finite. For $F\subseteq
\Bbb N$, let $\mu(F)=\sum_{i\in F}\mu(i)$. If $B\subseteq\Bbb N$,
then $\textrm{M}_B$ will denote the set of all convex means $\mu$
on $\Bbb N$ such that $\textrm{supp}(\mu)\subseteq B$. Let $\cal
F$ be a collection of finite subsets of $\Bbb N$. We denote
$\textrm{M}_B(\cal F,\epsilon)=\{\mu\in\textrm{M}_B:\forall F\in
\cal F~ \mu(F)<\epsilon\}$. Then

\begin{fct}[Pt\'{a}k's lemma] The two following are equivalent:

\begin{itemize}
             \item [(i)]  There exists a strictly increasing
             sequence $A_1\subset A_2\subset \cdots$ of finite
             subsets of $A$, and a sequence $F_n\in\cal F$ such
             that $F_n\subseteq A_n$ for all $n$.
             \item [(ii)] There exists an infinite subset $B$ of
             $A$ and an $\epsilon>0$ such that $\textrm{M}_B(\cal
             F,\epsilon)=\emptyset$.
\end{itemize}
\end{fct}
{\bf Proof.} See  \cite{Kot},   page 327. \qed

\begin{fct} \label{convex-interchangeable}  If a bounded subset $A$ of $C(X)$ has interchangeable double limits property, then so does
$\textrm{conv}(A)$.
\end{fct}
{\bf Proof.} Use Pt\'{a}k's lemma. The proof is a straightforward
translation of (4) in \cite{Kot}, page 328, for the topology of
pointwise convergence. \qed

\begin{cor}[The direction (iii) $\Rightarrow$ (i)] \label{(iii)->(i)}
If $A$ is countably weakly precompact, then $A$ is weakly
precompact.
\end{cor}
{\bf Proof.} This is a concequence of the above fact, the theorem
of  definability of types and Remark~\ref{rem}.
 Indeed, let $Y$ be any convex subset
of $C(X)$. If an arbitrary subset $A\subseteq Y$ is countably
precompact in $Y$, then $A$ is  countably precompact in $C(X)$.
By Fact \ref{convex-interchangeable}, $\textrm{conv}(A)$ is
 countably precompact in $C(X)$,
 so $\overline{\textrm{conv}(A)}$ , the closure of $\textrm{conv}(A)$ in $C(X)$, is
compact  (see Theorem~\ref{ES-pointwise}).
 Now by the theorem of definability of
types (Fact \ref{Definability} above), it is easy to verify that
$\overline{\textrm{conv}(A)}\subseteq Y$ and is the closure of
$\textrm{conv}(A)$ in $Y$ (see the paragraph after Fact
\ref{Definability}). Thus $\textrm{conv}(A)$ is precompact in $Y$,
so $\bar{A}$ is compact in $Y$. Therefore, by Remark~\ref{rem}
above, an arbitrary subset $A$ (of a normed space $Y$) is weakly
precompact if it is weakly countably precompact. \qed

\medskip

Fact \ref{convex-interchangeable} leads us to a new
characterization of local stability inside a model. Assume that
$\M$ is a structure, $\phi(x,y)$ a formula. Let
\begin{align*}
\textrm{conv}(\phi)=\{\psi(x,\bar{b}): \psi(x,\bar{b}) & \mbox{
is a convex combination of} \\ &  \mbox{  (at most finitely many)
formulas } \phi(x,b),b\in M \}.
\end{align*}

\begin{cor} Assume that $\phi(x,y)$ and $\M$ are as above. Then
the following are equivalent.
\begin{itemize}
             \item [(i)]  The formula $\phi$ is stable in $\M$.
             \item [(ii)] If $a_n\in\M$ and $\psi(x,\bar{b_n})\in\textrm{conv}(\phi)$ form two sequences we have
             $$\lim_n\lim_m\psi(a_m,\bar{b_n})=\lim_m\lim_n\psi(a_m,\bar{b_n}),$$
             whenever both limits exist.
\end{itemize}
\end{cor}

\medskip
By Mazur's lemma (see below), the theorem of  definability of
types is a consequence of Theorem \ref{ES-pointwise} above.
\begin{rmk}[Mazur Lemma] {\em If $(f_n)$ is a
bounded sequence of continuous functions on $X$ which converges
to a {\em continuous} function $f$, there exists a sequence
$g_n\in \textrm{conv}((f_k)_{k\geqslant n})$ which {\em
uniformly} converges to $f$. (Here $\textrm{conv}((h_k))$ denotes
the set of convex combinations of the $h_k$'s.) Therefore $f$ can
be written as a uniform limit of continuous functions of the form
$\frac{1}{n}\sum_{i<n}f_i$.  Standard proofs of Mazur's lemma use
the Hahn-Banach theorem and Lebesgue's Dominated Convergence
Theorem.
Pt\'{a}k's lemma gives a different proof of  Mazur's result
 (see \cite{Tod}, page 14). }
\end{rmk}

Thus, a formula $\phi(x,y)$ is stable in a model $M$ iff for every
type $p\in S_\phi(M)$ there is a sequence $\phi(a_n,y)$, $a_n\in
M$, (pointwise) convergening to a continuous function $\psi(y)$
such that $\phi(x,b)^p=\psi(b)$ for all $b\in M$. (Indeed, note
that $C(X)$ is a Fr\'{e}chet-Urysohn space (see Fact~
\ref{angel}).)



\noindent\hrulefill

\section{Model theory and Banach space theory}
In this short  section we continue the discussion of some the
topics raised above.

Recall that a Banach space $X$ is reflexive if a certain natural
isometry of $X$ into $X^{**}$ is onto. This mapping is
$\widehat{\ }:X\to X^{**}$ given by $\hat{x}(x^*)=x^*(x)$.

Now, we analyze the weak sequential compactness. Obviously, a
Banach space $X$ is weakly sequentially compact if the following
conditions hold:
\begin{itemize}
             \item [(a)] every bounded sequence $(x_n)$ of $X$ has a weak Cauchy subsequence,
             (i.e. there is $(y_n)\subseteq (x_n)$ such that for all $x^*\in X^*$ the sequence $(x^*(y_n))_{n=1}^\infty$ is a convergent sequence of reals, so
             $\hat{y}_n\to x^{**}$ weak* in $X^{**}$ for some $x^{**}\in X^{**}$),  and
             \item [(b)] every weak Cauchy sequence $(x_n)$ of $X$ has a weak limit (i.e. if $\hat{x}_n\to x^{**}$ weak* in $X^{**}$ then $x^{**}\in \widehat{X}$).
\end{itemize}

It is easy to check that the condition (a) corresponds to NIP and
the condition (b) corresponds to NSOP in model theory (see
below). In functional analysis, the condition (a) is called the
{\em weak* sequential compactness of $\widehat{X}$} (short
W*S-compactness), and the condition (b) is called the {\em weak
sequential completeness} (short WS-completetness). Clearly, a
weakly sequentially compact set is  weakly* sequentially compact,
but the converse fails. Indeed, the sequence $y_n=(
\underbrace{1,\ldots,1}_{n-times},0,\ldots)$ form a weakly Cauchy
sequence in $c_0$ without weak limit.

On the other hand, the Fr\'{e}chet-Urysohn property of the space
$C(X)$ corresponds to the definability of types in stable
theories: Let $\phi$ be a formula, $M$ a model of a stable theory
in continuous logic, and $S_\phi(M)$ the space of all complete
$\phi$-types on $M$ (see \cite{BBHU}). A type $p$ in $S_\phi(M)$
is a point in the closure of realized types in $M$, thus if
$a_\alpha\in M$ and $tp(a_\alpha/M)\to p$, then there exists a
continuous function $\psi$ such that $\phi(a_\alpha,y)\to\psi(y)$
pointwise ($\psi$ is continuous because the theory is stable, see
\cite{K3}). Now, by the Fr\'{e}chet-Urysohn property of
$C(S_\phi(M))$, there is a sequence $(a_n)\subseteq M$ such that
$\lim_n\phi(a_n,y)=\psi(y)$, i.e. $p$ is definable by $\psi$.

A standard fact in functional analysis is that a Banach space $X$
is reflexive if and only if $B_X=\{x\in X:\|x\|\leqslant 1\}$ is
weakly compact. Thus by the Eberlain-\v{S}mulian theorem, $X$ is
reflexive if and only if the conditions (a) and (b) above hold
for $A=B_X$. This and the above observations show that stability
in model theory corresponds to reflexivity in functional analysis.
Thus, one can say that `first order logic is angelic.' To
summarize:

$$\begin{tabular}{|c|c|}
  \hline
  Logic & Analysis \\ \hline  \hline
  \ \ \ \ Definability of types \ \ \ \  & Fr\'{e}chet-Urysohn property \\ \hline
  NIP & W*S-compactness  \\ \hline
  NSOP & WS-completeness \\ \hline
  Shelah's theorem & Eberlein-\v{S}mulian theorem \\ \hline
\end{tabular}$$

\bigskip
In the next subsection we study these connections more exactly.

\subsection{Type space}
We assume that the reader is familiar with continuous logic (see
\cite{BBHU}). Of course, we study  real-valued formulas instead of
$[0,1]$-valued formulas. One can assign bounds to formulas and
retain compactness theorem in a local way again.

Suppose that $L$ is an arbitrary language.  Let $M$ be an
$L$-structure, $A\subseteq M$ and $T_A=Th({M}, a)_{a\in A}$. Let
$p(x)$ be a set of $L(A)$-statements in free variable $x$. We
shall say that $p(x)$ is a {\em type  over} $A$ if $p(x)\cup T_A$
is satisfiable. A {\em complete type over} $A$ is a maximal type
over $A$. The collection of all such types over $A$ is denoted by
$S^{M}(A)$, or simply by $S(A)$ if the context makes the theory
$T_A$ clear. The {\em type of $a$ in $M$ over $A$}, denoted by
$\text{tp}^{M}(a/A)$, is the set of all $L(A)$-statements
satisfied in $M$ by $a$. If $\phi(x,y)$ is a formula, a {\em
$\phi$-type} over $A$ is a maximal consistent set of formulas of
the form $\phi(x,a)\geqslant r$, for $a\in A$ and
$r\in\mathbb{R}$. The set of $\phi$-types over $A$ is denoted by
$S_\phi(A)$.

We now give a characterization of complete types in terms of
functional analysis. Let $\mathcal{L}_A$ be the family of all
interpretations $\phi^{M}$ in $M$ where $\phi$ is an
$L(A)$-formula with a free variable $x$. Then $\mathcal{L}_A$ is
an Archimedean Riesz space of measurable functions on $M$ (see
\cite{Fremlin3}). Let $\sigma_A({M})$ be the set of Riesz
homomorphisms $I: {\mathcal L}_A\to \mathbb{R}$ such that
$I(\textbf{1}) = 1$. The set $\sigma_A({M})$ is   called the {\em
spectrum} of $T_A$. Note that $\sigma_A({M})$ is a weak* compact
subset of $\mathcal{L}_A^*$. The next proposition shows that a
complete type can be coded by a Riesz homomorphism and gives a
characterization of complete types. In fact, by Kakutani
representation theorem, the map $S^{M}(A)\to\sigma_A({M})$,
defined by $p\mapsto I_p$ where $I_p(\phi^M)=r$ if $\phi(x) = r$
is in $p$, is a bijection.

\begin{pro} \label{key}
Suppose that $M$, $A$ and $T_A$ are as above.
\begin{itemize}
             \item [{\em (i)}] The map $S^{M}(A)\to\sigma_A({M})$  defined by $p\mapsto I_p$ is bijective.
             \item [{\em (ii)}] $p\in S^{M}(A)$ if and only if there is an elementary
extension $N$ of $M$ and $a\in N$ such that
$p=\text{tp}^{N}(a/A)$.
\end{itemize}
\end{pro}

We equip $S^{M}(A)=\sigma_A({M})$  with the related topology
induced from $\mathcal{L}_A^*$. Therefore, $S^{M}(A)$ is a compact
and Hausdorff space. For any complete type $p$ and formula
$\phi$, we let $\phi(p)=I_p(\phi^{M})$. It is easy to verify that
the topology on $S^{M}(A)$ is the weakest topology in which all
the functions $p\mapsto \phi(p)$ are continuous. This topology
sometimes called the {\em logic topology}. The same things are
true for $S_\phi(A)$.

It is easy to check that for each $I\in\sigma_M({\M})$, $\|I\|=1$
and $I$ is a positive and multiplicative, i.e. $I(f\times
g)=I(f)\times I(g)$ for all $f,g\in\mathcal{L}_M$. (Recall that
$I$ is positive if $I(f)\geqslant 0$ for all $f\geqslant 0$.)

Also, $\sigma_M({\M})$ is the set of extreme points of the state
space $\mathcal{S}=\{I\in \mathcal{L}_M^*:I(1)=1\mbox{ and $I$ is
positive }\}$. $\mathcal{S}$ also called the space of Keisler
types. (By the Krein-Milman theorem,
$\mathcal{S}=\overline{conv}^{w}(\sigma_M({\M}))$ , i.e. it is
the (weak) closure of convex hull of its extreme points.) Then
$\sigma_M({\M})\subset \mathcal{S}\subset \mathcal{S}^1=\{I\in
\mathcal{L}_M^*:\|I\|= 1\}\subset B_{\mathcal{L}_M^*}=\{I\in
\mathcal{L}_M^*:\|I\|\leqslant 1\}$.

Since every member of $\mathcal{L}_M^*$ is expressible as the
difference of two positive linear functionals, $\sigma_M({\M})$
determines the set $B_{\mathcal{L}_M^*}$ (and hence the Banach
space $\mathcal{L}_M^*$).  Thus, space of types can be equipped
with a natural  norm space structure, and we can study this
Banach space (i.e. $\mathcal{L}_M^*$) instead of the types space.
The weak* topology of this Banach space on the space of types is
the logic topology, and we have a natural linear structure on the
space of types, i.e. for all types $p,q$, the addition $p+q$ is
well defined, also $rp$ is well defined for each real number $r$.
(Indeed, $p+q:=I_p+I_q$ and $rp:=rI_p$.) Of course, $p+q$ is not
necessary a (classical) type, but it is easier to study the
Banach space determined by types.

\subsection{Banach space for a formula}
Let $\M$ be an $L$-structure, $\phi(x,y):M\times M\to \Bbb R$ a
formula (we identify formulas with real-valued functions defined
on models).

Let $S_\phi(M)$ be the space of complete $\phi$-types over $M$ and
set $A=\{\phi(x,a),-\phi(x,a)\in C(S_\phi(M)): a\in M\}$. The
(closed) convex hull of $A$, denoted by ($\overline{conv}(A)$)
$conv(A)$, is the intersection of all (closed) convex sets that
contain $A$. $\overline{conv}(A)$ is convex and closed, and
$\|f\|\leqslant \|\phi\|$ for all $f\in\overline{conv}(A)$. So,
by normalizing we can assume that $\|f\|\leqslant 1$ for all
$f\in\overline{conv}(A)$. We claim that $B=\overline{conv}(A)$ is
the unit ball of a Banach space. Set
$V=\bigcup_{\lambda>0}\lambda B$. It is easy to verify that $V$
is a Banach space with the normalized norm and $B$ is its unit
ball. This space will be called the {\em space of linear
$\phi$-definable relations}. One can give an explicit description
of it: $$V=\Big(~\overline{\big\{\sum_1^n r_i\phi(x,a_i):a_i\in M,
r_i\in{\Bbb R},n\in{\Bbb N}\big\}}~;~~\||\cdot\|| ~ \Big)$$ where
$ \||\cdot\||$ is the normalized norm.

\medskip
Note that $V$ is a subspace of $C(S_\phi(M))$. Recall that for an
infinite compact Hausdorff $X$, the space  $C(X)$ is not
reflexive, nor is it weakly complete.  So, if $V$ is a lattice
(or algebra), then it is not reflexive, nor is it weakly complete
(since, in this case,  $V$ is isomorphic to $C(X)$ for some
compact Hausdorff space $X$).

\subsection{Stability and reflexivity}
 A formula $\phi:M\times M\to {\Bbb R}$ has the order property  if there
exist sequences $(a_m)$ and $(b_n)$ in $M$ such that
$$\lim_m\lim_n\phi(a_m,b_n)\neq \lim_n\lim_m\phi(a_m,b_n)$$
We say that  $\phi$ has the double limit property (DLP) if it has
not has the order property.

\begin{dfn} We say that $\phi(x,y)$ is unstable if either $\phi$
or $-\phi$ has the order property. We call $\phi$ stable if
$\phi$ is not unstable.
\end{dfn}

The following is a well known result in functional analysis:

\begin{fct}
A Banach space $X$ is reflexive iff
 its unite ball is weakly compact.
\end{fct}

If $\phi$ is stable then the set $A$ (see above) is weakly
precompact in $C(S_\phi(M))$ by Grothendiek's criterion. (Note
that the collection of types realized in $\M$ is dense in
$S_\phi(M)$.) By Pt\'{a}k's lemma, the convex hull of $A$,
$conv(A)$, is weakly precompact, and therefore the closed convex
hull  of $A$, i.e. $B$, is so. (Note that for convex sets, weakly
closed $=$ uniform closed.)

\begin{cor} 
Assume that $\phi(x,y)$, $\M$, $B$ and $V$ are as above. Then the
following are equivalent:

(i) $\phi$ is stable on $\M$.

(ii)  $B$  is weakly compact.

(iii)   The Banach space $V$ is reflexive.
\end{cor}

Recall that for an infinite compact Hausdorff space $X$, the
space $C(X)$ is not reflexive.

\subsection{NIP and Rosenthal Banach spaces}
We say that a formula $\phi(x,y)$ is NIP on a model $\M$ if every
sequence of the  set $A=\{\phi(x,a),-\phi(x,a)\in C(S_\phi(M)):
a\in M\}$ is not independent in the sense of Definition
\ref{indep}.

\begin{dfn}[\cite{GM}, 2.10] {\em A Banach space $X$ is said to be {\em Rosenthal} if it
does not contain an isomorphic copy of $\ell_1$.}
\end{dfn}

\begin{fct}[H.P. Rosenthal]
For a Banach space $X$ the following either $X$ is Rosenthal or
it has a sequence with no weak Cauchy subsequence.
\end{fct}

If $\phi(x,y)$ is NIP on $\M$ then the set $A$ is not an
independent family, and $conv(A)$ is nor (see \cite{BFT}, page
878). (Note that the collection $M_0$ of types realized in $\M$
is dense in $S_\phi(M)$ and $A$ is not independent iff the family
$\{f|_{M_0}:f\in A\}$ is not independent (see \cite{GM}, Lemma
7.9.) So, its (uniform) closed convex hull, $\overline{conv}(A)$,
is not independent. Now, the above fact and Lemma~\ref{lemma-1}
imply that:

\begin{cor}
Assume that $\phi(x,y)$, $\M$, and $V$ are as above. Then the
following are equivalent:

(i) $\phi$ is NIP on $\M$.

(ii) $V$ is Rosenthal Banach space.

(iii)  Every bounded sequence of $V$ has a
 weak Cauchy subsequence.
\end{cor}

\subsection{NSOP and weak sequential completeness} Let ${\M}(={\mathcal U})$ be
a monster model (of theory $T$) and $\phi(x,y)$ a formula.

\begin{cor}
Assume that $\phi(x,y)$, $\M$, and $V$ are as above. Then the
following are equivalent:

(i) $\phi$ is NSOP (on $\M$).

(ii)  Every weak Cauchy sequence of $V$ has a weak limit (in $V$).

\end{cor}
{\bf Proof.}  Use the above corollaries and the
Eberlein-\v{S}mulian Theorem. \qed

\medskip

Note that for a compact Hausdorff space $X$, the space $C(X)$
contains an isomorphic copy of $c_0$, and so $C(X)$ is not weakly
sequentially complete (see \cite{AK}, Proposition 4.3.11).

\medskip
In \cite{Iov} Jose Iovino pointed out the correspondence between
stability and reflexivity. He showed that a formula $\phi(x,y)$
is stable iff $\phi$ is the pairing map on the unite ball of
$E\times E^*$, where $E$ is a reflexive Banach space. In this
paper, we gave a `concrete and explicit' description of the
Banach space $V$, such that it is reflexive iff $\phi$ is stable.
This space is uniquely determined by $\phi$ and the formula
$\phi$ is completely coded by $V$. The value of $\phi$ is exactly
determined by the evaluation map $\langle
\cdot,\cdot\rangle:V\times V^*\to\Bbb R$ defined by $\langle
f,I\rangle=I(f)$.

\bigskip
At the end of this paper - of not the story - we note that the
basic ideas came from model theory, in the other word, techniques
from one field became useful in the other. One might therefore
hope to obtain other connections between stability theory and
Banach space theory. We will continue this way in a future work.

\noindent\hrulefill

\appendix

\section{Stability and Eberlein compacta}

{\bf Proof of Fact~\ref{Criterion}.} (i)~$\Rightarrow$~(ii):
Suppose that $f_n\in A$ and $x_n\in X$ form two sequences and the
limits $\lim_n \lim_m f_n(x_m)$ and $\lim_m \lim_n f_n(x_m)$
exist. Let $f$ in $C(X)$ and $x$ in $X$ be cluster points of
$\{f_n\}$ and $\{x_m\}$. Thus,
\begin{align*} \lim_n\lim_m f_n(x_m)  =
\lim_nf_{n}(x) &= f(x)=\lim_mf(x_m)=\lim_m\lim_nf_n(x_m).
\end{align*}

(ii)~$\Rightarrow$~(i): Since $A$ is bounded (i.e. there is an
$r$ such that $|f|\leqslant r$ for all $f\in A$) and by
Tychonoff's theorem $[-r,r]^X$ is compact, it suffices to show
that $\overline{A}\subseteq C(X)$. Let $f\in \overline{A}$.
Suppose that $f$ is not continuous at a point $x$ in $X$. Then
there is a neighborhood $U$ of $f(x)$ such that each neighborhood
of $x$ contains a point $y$ of $X$ with $f(y)$ not belonging to
$U$. Take any $f_1$ in $A$; then there is an $x_1$ in $X$ such
that $|f_1(x)-f_1(x_1)|<1$ and $f(x_1)\notin U$. Take $f_2$ in
$A$ so that $|f_2(x_1)-f(x_1)|<1$ and $|f_2(x)-f(x)|<1$. Now
choose $x_2$ in $X$ such that $|f_i(x)-f_i(x_2)|< 1/2$ ($i=1,2$)
and $f(x_2)\notin U$. Then take $f_3$ in $A$ so that
$|f(x_j)-f_3(x_j)|< 1/2$ and $|f(x)-f_3(x)|< 1/2$. Proceeding in
this way, one obtains sequences $\{f_n\}_n$ and $\{x_m\}_m$ in
$A$ and $X$ such that, for each $n$, $|f_i(x)-f_i(x_n)|<1/n$
($i=1, 2,\ldots,n$), $|f(x_j)-f_{n+1}(x_j)|<1/n$
($j=1,2,\ldots,n$), $|f(x)-f_{n+1}(x)|< 1/n$, and $f(x_n)\notin
U$. Then $\lim_n \lim_m f_n(x_m) =\lim_nf_n(x) = f(x)$, and
$\lim_nf_n(x_m)=f(x_m)\notin U$. Since it is possible to take a
subsequence of $\{x_m\}_m$ so that the corresponding subsequence
of $\{f(x_m)\}_m$ converges to a point outside of $U$, the
assumption that $f$ is not continuous contradicts the iterated
limit condition of (ii). \qed

\noindent\hrulefill

\section{Angelicity of $C(X)$}
\begin{dfn}[Fremlin] {\em A regular Hausdorff space space is {\em angelic} if (i) every
countably precompact set is precompact, (ii) the closure of a
precompact set is precisely the set of limits of its sequences.}
\end{dfn}

For angelicity of $C(X)$ (where $X$ is a compact Hausdorff space)
it suffices to so prove that:

\begin{fct}[\cite{Fremlin4}, 462B]  \label{angel}
Assume that $X$ is compact and $A\subseteq C(X)$ is bounded and
precompact. If $g\in\bar{A}$, there is a sequence $f_n\in A$ such
that $\lim_n f_n=g$.
\end{fct}
{\bf Proof.} For $g\in\overline{A}$ we construct  countable sets
$D\subseteq X$, $B\subseteq A$ such that

 (1) whenever $I\subseteq B\cup\{g\}$ is finite, $\epsilon>0$
and $x\in X$, there is a $y\in D$ such that
$|f(y)-f(x)|\le\epsilon$ for every $f\in I$;

 (2) whenever $J\subseteq D$ is finite and $\epsilon>0$ there
is an $f\in B$ such that $|f(x)-g(x)|\le\epsilon$ for every $x\in
J$.

  For any finite set $I\subseteq\Bbb R^X$, the set
$Q_{I}=\{\{ f(x)\}_{f\in I}:x\in X\}$ is a subset of the
separable metrizable space $\Bbb R^I$, so is itself separable,
and there is a countable dense set $D_{I}\subseteq X$ such that
$Q'_{I}=\{\{ f(x)\}_{f\in I}:x\in D_{I}\}$ is dense in $Q_{I}$.
Similarly, because $g\in\overline{A}$, we can choose for any
finite set $J\subseteq X$ a sequence $\{f_{Ji}\}_i$ in $A$ such
that $\lim_if_{Ji}(x)=g(x)$ for every $x\in J$.

Now construct $\{D_m\}_m$, $\{B_m\}_m$ inductively by setting
$D_m=\bigcup\{D_{I}: I\subseteq\{g\}\cup\bigcup_{i<m}B_i$ is
finite$\}$, $B_m=\{f_{Jk}:k\in\Bbb N,\,J\subseteq\bigcup_{i<m}D_i$
is finite$\}$.
 At the end of the induction, set $D=\bigcup_{m\in\Bbb
N}D_m$ and $B=\bigcup_{m\in\Bbb N}B_m$.   Since the construction
clearly ensures that $\{D_m\}_m$ and $\{B_m\}_m$ are
non-decreasing sequences of countable sets, $D$ and $B$ are
countable, and we shall have $D_{I}\subseteq D$ and $I\subseteq
B\cup\{g\}$ is finite, while $f_{Ji}\in B$ whenever $i\in\Bbb N$
and $J\subseteq D$ is finite.   Thus we have suitable sets $D$
and $B$.

By (2), there must be a sequence $\{f_i\}_i$ in $B$ such that
$g(x)=\lim_if_i(x)$ for every $x\in D$.   In fact
$g(y)=\lim_if_i(y)$ for every $y\in X$. Otherwise, there is an
$\epsilon>0$ such that $J=\{i:|g(y)-f_i(y)|\ge\epsilon\}$ is
infinite.    For each $m\in\Bbb N$, $I_m=\{f_i:i\le m\}$ is a
finite subset of $B$, so by (1) there is an $x_m\in D_{I_m}$ such
that $|f(x_m)-f(y)|\le 2^{-m}$ for every $f\in I_m\cup\{g\}$. Then






$$\lim_i\lim_mf_i(x_m)=\lim_if_i(y)\neq g(y)=\lim_mg(x_m)=\lim_m\lim_if_i(x_m).$$

But this is impossible, by Fact~\ref{Criterion}. So $g=\lim_if_i$.  
\qed

\medskip

Now, since any subspace of an angelic space is angelic (see
\cite{Fremlin4}, 462C(a)), by Remark~\ref{rem} above, every
Banach space (or even normed space) is angelic. This leads to
another proof of  the direction (iii)~$\Rightarrow$~(i) of the
Eberlain-\v{S}mulian (see Corollary~\ref{(iii)->(i)}).


\noindent\hrulefill


\end{document}